\documentclass[11pt]{amsart}
\usepackage[T1]{fontenc}
\usepackage{lmodern}
\usepackage{xcolor}
\usepackage{amsmath,amssymb,amsthm,mathtools}
\usepackage[margin=1in]{geometry}
\usepackage{microtype}
\usepackage{enumitem}
\usepackage[hidelinks]{hyperref}

\newtheorem{theorem}{Theorem}[section]
\newtheorem{proposition}[theorem]{Proposition}
\newtheorem{lemma}[theorem]{Lemma}

\theoremstyle{remark}
\newtheorem{remark}[theorem]{Remark}
\numberwithin{equation}{section}
\newcommand{\C}{\mathbb C}
\newcommand{\PP}{\mathbb P}
\newcommand{\N}{\mathbb N}
\newcommand{\supp}{\operatorname{supp}}
\newcommand{\diam}{\operatorname{diam}}
\newcommand{\Span}{\operatorname{span}_{\C}}
\newcommand{\Dim}{\dim_{\mathrm H}}
\newcommand{\BD}{\mathcal D_{\mathrm P}}

\title[Base-dominant rigidity]{History independence and base-dominant rigidity\texorpdfstring{\\}{ }for polynomial skew products}
\author{Yan Mary He}
\dedicatory{To Peter Shalen on his 80th birthday, with gratitude and admiration.}
\address{Department of Mathematics\\
University of Oklahoma\\
Norman, OK 73019, USA}
\email{he@ou.edu}
\date{\today}

\begin{document}

\begin{abstract}
We study history independence of fast directions and its relation to the
Binder--DeMarco dimension formula for regular polynomial skew products
of degree $d\ge2$.
Assuming uniform expansion on the support of the measure of maximal
entropy and a base Lyapunov exponent larger than the vertical exponent,
we prove that history independence is equivalent to conjugacy to a
polynomial product by an affine shear map.
Under history independence, the measure of maximal entropy satisfies the
Binder--DeMarco dimension formula.
If, in addition, the splitting is dominated, the dimension formula
also implies history independence.
The proof first upgrades the measurable fast direction to a holomorphic
line field and then uses its invariance on periodic fibers to force
the affine product structure.
The converse dimension implication combines an entropy comparison argument with
rigidity of local symmetries of periodic fiber Julia sets.
\end{abstract}

\maketitle

\section{Introduction}
Let $f$ be a holomorphic endomorphism of $\PP^k=\PP^k(\C)$
of algebraic degree $d\ge2$, and let $\mu_f$ be its unique measure
of maximal entropy \cite{Lyubich1983,FreireLopesMane1983,BriendDuval2001}.
Write $\chi_1\ge\cdots\ge\chi_k>0$ for the Lyapunov exponents of $\mu_f$,
counted with complex multiplicity \cite{BriendDuval1999}.
In dimension one, the Ma\~n\'e--Manning formula
\cite{Mane1988,Manning1984} gives $\Dim\mu_f=\log d/\chi_1$.
Binder--DeMarco \cite{BD} conjectured that the Hausdorff dimension $\Dim\mu_f$ of $\mu_f$ satisfies
$\Dim\mu_f=\sum_{j=1}^k\log d/\chi_j$.

In dimension two, when the Lyapunov exponents are distinct,
the fast Oseledec line is defined along almost every inverse
history and may vary between histories ending at the same point.
Bianchi and the author \cite{BH} proved that, under hypothesis
\textup{(H)} of \cite{BH} and with $\chi_1>\chi_2>\log d$,
history dependence implies
$\Dim\mu_f>\log d/\chi_1+\log d/\chi_2$.
They constructed a non-empty open set of endomorphisms of $\PP^2$
satisfying these conditions, thereby disproving the Binder--DeMarco
conjecture. Adjoining independent quadratic polynomial factors to their  examples gives
counterexamples in every dimension $k\ge3$.

In this paper, we study the relation between history independence
and the Binder--DeMarco formula for regular polynomial skew products.
Let $F(z,w)=(p(z),q(z,w))$ be a degree $d\ge2$ polynomial skew product
extending holomorphically to $\PP^2$. Write $J_p$ for the Julia set of $p$.
Let $\mu:=\mu_F$ be the unique measure of maximal entropy, and set
$K:=\supp\mu$. We refer to Jonsson \cite{Jonsson} for background on polynomial skew products and to Astorg--Bianchi \cite{ABHigher,AB} for
hyperbolicity and bifurcations in holomorphic families. We write
$$
 \lambda_b:=\int\log|p'(z)|\,d\mu(z,w),\qquad
 \lambda_v:=\int\log|q_w(z,w)|\,d\mu(z,w)
$$
for the base and vertical Lyapunov exponents, and set
$\BD(F):=\log d/\lambda_b+\log d/\lambda_v$.

The fiber-dominant case $\lambda_v>\lambda_b$ is straightforward:
the fast direction is vertical and hence history-independent.
The formula $\Dim\mu_F=\BD(F)$ was proved in
\cite[Proposition 5.3]{BH} for an explicit family satisfying
\textup{(H)}.
We observe that the Ledrappier--Young dimension argument
used there yields the formula throughout the fiber-dominant
case, without uniform expansion or domination; we give the
details in Appendix \ref{sec:fibered}.
The general formula was also proved independently by Genes \cite{Genes} using a different method.
We therefore focus on the base-dominant case $\lambda_b>\lambda_v$.
It turns out that there is an unexpected {\it rigidity} phenomenon: under uniform
expansion, measurable history independence forces conjugacy to a
polynomial product by an affine shear map.

Write $\hat K$ for the natural extension of $F|_K$, with points
$\hat x=(x_0,x_{-1},\ldots)$ satisfying $F(x_{-j})=x_{-j+1}$,
write $\hat F$ for the lifted map, and set $\pi_0(\hat x):=x_0$.
Let $\hat\mu$ be the lift to the natural extension of $\mu$.

\smallskip
\noindent\textbf{Expanding setting (E$_0$).}
The map $F$ is uniformly expanding on $K$: there are $c>0$ and $\rho>1$
such that $\|DF_x^m v\|\ge c\rho^m\|v\|$ for every $x\in K$,
$v\in T_x\C^2$, and $m\ge1$. Moreover, $\lambda_b>\lambda_v$.

Let $E^f(\hat x)$ denote the fast Oseledec line.
We call the fast direction \emph{history-independent} if there
is a measurable line field $\ell$ on $K$ such that
$E^f(\hat x)=\ell(x_0)$ for $\hat\mu$-almost every $\hat x$.

\begin{theorem}\label{thm:rigidity}
Suppose that $F$ satisfies \textup{(E$_0$)}.
Then the fast direction is history-independent if and only if there are a
constant $a\in\C$ and a degree $d$ polynomial $Q$ such that
\begin{equation}\label{eq:product-form}
 q(z,w)=a\,p(z)+Q(w-az).
\end{equation}
Equivalently, the affine shear map $H_a(z,w):=(z,w-az)$ satisfies
$H_a\circ F\circ H_a^{-1}(z,v)=(p(z),Q(v))$.
In this case $E^f(\hat x)=\Span(1,a)$
for $\hat\mu$-almost every $\hat x$.
\end{theorem}

Under history independence, Theorem \ref{thm:rigidity} and the Ma\~n\'e--Manning formula
\cite{Mane1988, Manning1984} imply that
$\Dim\mu_F=\BD(F)$.
To characterize equality, we also assume a dominated splitting.

\smallskip
\noindent\textbf{Expanding dominated setting (E).}
In addition to \textup{(E$_0$)}, the natural extension admits a continuous
invariant splitting $T_{x_0}\C^2=E^f(\hat x)\oplus E^v(x_0)$,
where $E^v:=\{0\}\times\C$, and constants $C>0$, $0<\tau<1$ such that
\begin{equation}\label{eq:domination}
 \frac{\|DF^m|_{E^v(x_0)}\|}{\|DF^m|_{E^f(\hat x)}\|}
 \le C\tau^m, \qquad \hat x\in\hat K, m\ge1.
\end{equation}
Our next result characterizes Binder--DeMarco equality under the
expanding dominated setting.

\begin{theorem}
\label{thm:BD}
Suppose that $F$ satisfies \textup{(E)}.
Then $\mu_F$ is exact-dimensional, and the following are equivalent:
\begin{enumerate}[label=\textup{(\roman*)},leftmargin=2em]
\item $\Dim\mu_F=\BD(F)$;
\item the fast direction is history-independent;
\item $F$ has the product form \eqref{eq:product-form}.
\end{enumerate}
\end{theorem}

For regular polynomial skew products, both Lyapunov exponents are
greater than or equal to $\log d$; see \cite{Jonsson} and \cite[Section 3.1]{AB}.
Under \textup{(E)}, history dependence is equivalent to
$\Dim\mu_F>\BD(F)$, even when $\lambda_v=\log d$.
Both equality and strict inequality occur at this boundary value $\lambda_v=\log d$;
see Remark \ref{rem:boundary}.

\medskip

The first step in proving Theorem \ref{thm:rigidity} is to show that
the measurable fast slope is holomorphic near $K$.
Invariance of the fast line field under $DF$ gives a cohomological
equation for its slope.
We construct local holomorphic solutions by convergent series along
inverse histories and use history independence to glue the local solutions near $K$.
We prove an extension lemma which shows that, under suitable degree and multiplier conditions,
holomorphic solutions of a polynomial cohomological equation near a
hyperbolic Julia set extend to entire affine functions.
After centering the fiber polynomials, we apply this lemma on
non-resonant periodic fibers to force the fast slope to vanish,
yielding the product structure.

The holomorphic fast line field locally defines an invariant foliation
near $K$. Favre--Pereira \cite{FP} classify non-invertible dominant
rational maps preserving a holomorphic foliation on a projective surface.
Tapiero \cite{Tapiero} obtains invariant foliations under
pluripotential product hypotheses. However, here we derive the foliation from \emph{measurable} history independence
and prove the affine product form directly.

For Theorem \ref{thm:BD}, the
Ledrappier--Young formula for endomorphisms reduces Binder--DeMarco
equality to $h_1=\log d$, where $h_1$ is the partial entropy along fast leaves.
Adapting the reference measure entropy argument of
Buzzi--Fisher--Tahzibi \cite{BFT}, we show that $h_1=\log d$ forces
every local fast plaque, restricted to points projecting into $J_p$,
to lie in $K$.
We compare fast holonomies associated with different inverse histories
from a periodic fiber to the same nearby fiber over $J_p$.
Composing one holonomy with the inverse of the other gives a local
symmetry of the periodic fiber Julia set.
Levin's rigidity theorem for these symmetries, together
with the circle case, forces history independence.

\subsection*{Organization of the paper}
In Section 2, we prove the extension lemma which will be used to prove
Theorem \ref{thm:rigidity} in Section 3.
In Section \ref{sec:base-dominant}, we prove Theorem \ref{thm:BD}.
In Appendix \ref{sec:fibered}, we treat the fiber-dominant case
using the Ledrappier--Young argument of \cite[Section 5.3]{BH}.

\subsection*{Acknowledgements}
The author is supported by the National Science Foundation grant
DMS-2554272 and by the Simons Foundation Travel Support for
Mathematicians program.

\subsection*{AI Declaration}
ChatGPT 6 was used to assist with numerical investigations, literature searches, the discussion
and checking of mathematical arguments, editorial revision and drafting the appendix whose arguments are similar to those in \cite{BH}.

\section{An extension lemma}
The goal of this section is to prove Lemma \ref{lem:extension}, which states that
a holomorphic solution of a polynomial cohomological equation near a
hyperbolic Julia set extends to an entire affine function.
For a polynomial $P$, denote its Julia set by $J_P$.
Write $D(a,r):=\{w\in\C:|w-a|<r\}$.
A holomorphic germ along a compact set is represented
by a holomorphic function on an open neighborhood of that set; two
representatives are identified after restriction to a smaller neighborhood.

\begin{lemma}\label{lem:extension}
Let $P(w)=a_Dw^D+O(w^{D-2})$ be a hyperbolic polynomial, where
$a_D\ne0$ and $D\ge2$. Suppose that a holomorphic germ $u$ along
$J_P$, a constant $\Lambda\in\C$ with $|\Lambda|>1$, and a polynomial
$R$ of degree at most $D-2$ satisfy
\begin{equation}\label{eq:cohomology}
 \Lambda u(P(w))=P'(w)u(w)+R(w)
 \quad\text{as germs along }J_P.
\end{equation}
Then $u$ extends uniquely to an entire function, and there is a
constant $A\in\C$ such that
$$
 u(w)=Aw,\qquad (\Lambda-D)A=0,\qquad
 R(w)=A\bigl(DP(w)-wP'(w)\bigr).
$$
In particular, $\Lambda\ne D$ implies $u=R=0$.
\end{lemma}

\begin{proof}
\emph{Step 1: Exhaustion away from attracting cycles.}
Let $\mathcal A$ be the finite set of finite attracting periodic
points of $P$. Since $P$ is hyperbolic, every point of its Fatou
set tends to $\mathcal A$ or to infinity; see
\cite[Section 3.4]{McMullen}. Choose a bounded open neighborhood
$\Omega$ of $J_P$ such that
$\overline\Omega\cap\mathcal A=\emptyset$,
$P^{-1}(\Omega)\Subset\Omega$, and $u$ is holomorphic on
$\Omega$ with \eqref{eq:cohomology} valid on $P^{-1}(\Omega)$.
Such neighborhoods can be chosen arbitrarily small. Indeed, let
$V\subset\PP^1$ be an open neighborhood of
$\mathcal A\cup\{\infty\}$ such that $P(\overline V)\subset V$
and $\overline V\cap J_P=\emptyset$.
As $N\to\infty$, the sets
$\PP^1\setminus P^{-N}(\overline V)$ form arbitrarily small
neighborhoods of $J_P$ with the required properties.

For $n\ge1$, set
$
A_n:=\{v\in\C:P^{-n}(v)\subset\Omega\}.
$
The set $A_n$ is open since its complement is the image of the closed set
$\C\setminus\Omega$ under the proper map $P^n$. Since
$P^{-1}(\Omega)\subset\Omega$, we have
$\Omega\subset A_n\subset A_{n+1}\subset\C\setminus\mathcal A$.
The last inclusion follows by taking inverse images along each
attracting cycle.

Every compact set $E\Subset\C\setminus\mathcal A$ is contained
in $A_n$ for all sufficiently large $n$. Indeed, otherwise there would be
$n_j\to\infty$ and $w_j\notin\Omega$ with $P^{n_j}(w_j)\in E$.
Choose $r_*>0$ so that $E\subset D(0,r_*)$ and
$|P(w)|>2|w|$ for $|w|\ge r_*$.
Since $P^{n_j}(w_j)\in E$, we have $|w_j|<r_*$.
Thus all $w_j$ lie in the compact subset
$L:=\overline{D(0,r_*)}\setminus\Omega$ of the Fatou set.
Choose an open neighborhood $W\subset\PP^1$ of
$\mathcal A\cup\{\infty\}$ with $P(W)\subset W$ and
$W\cap E=\emptyset$. Since $L$ lies in the Fatou set and $W$
contains all attracting cycles and infinity, every point of $L$
eventually enters $W$. Compactness of $L$, together with $P(W)\subset W$,
then gives an integer $N$ such that $P^n(L)\subset W$ for every $n\ge N$.
This contradicts $P^{n_j}(w_j)\in E$.

We now construct a nested exhaustion of $\C\setminus\mathcal A$
by connected domains containing $\Omega$.
Choose a compact connected set
$E_0\subset\C\setminus\mathcal A$ containing $\overline\Omega$.
For all sufficiently large $n$, let $\Omega_n$ be the connected
component of $A_n$ containing $E_0$. These components are nested
and exhaust $\C\setminus\mathcal A$. To see the exhaustion,
join any point of $\C\setminus\mathcal A$ to $E_0$ by a compact
path avoiding $\mathcal A$. The union of the path and $E_0$ is a
compact subset of $\C\setminus\mathcal A$, so it lies in $A_n$
for all sufficiently large $n$.

\smallskip
\noindent\emph{Step 2: Extension to $\C\setminus\mathcal A$.}
Since $P^{-1}(\Omega)\subset\Omega$, every $w\in P^{-n}(\Omega)$
satisfies $P^j(w)\in\Omega$ for $0\le j\le n$.
We may therefore iterate \eqref{eq:cohomology} to obtain polynomials
$R_n$ such that
\begin{equation}\label{eq:iterated-cohomology}
 \Lambda^n u(P^n(w))=(P^n)'(w)u(w)+R_n(w),
 \qquad w\in P^{-n}(\Omega).
\end{equation}

For $v\in A_n$, define
$$
 u_n(v):=\frac{1}{D^n\Lambda^n}
 \sum_{P^n(w)=v}\bigl((P^n)'(w)u(w)+R_n(w)\bigr),
$$
where inverse images are counted with multiplicity. The function $u_n$
is holomorphic on $A_n$, including at critical values of $P^n$.
To see this, fix $v_0\in A_n$
and take small disjoint circles in $\Omega$ enclosing all distinct
roots of $P^n(w)-v_0$. The union $\Gamma$ of these circles encloses all roots for
$v$ near $v_0$. Setting $H_n(w):=(P^n)'(w)u(w)+R_n(w)$, we have
$$
 u_n(v)=\frac{1}{2\pi iD^n\Lambda^n}\int_\Gamma
 H_n(w)\frac{(P^n)'(w)}{P^n(w)-v}\,dw,
$$
which is holomorphic in $v$. The factor $(P^n)'(w)/(P^n(w)-v)$
has residue which is equal to the multiplicity of each root.

Equation \eqref{eq:iterated-cohomology} gives $u_n=u_{n+1}=u$ on $\Omega$.
The identity theorem therefore gives
$u_{n+1}|_{\Omega_n}=u_n|_{\Omega_n}$.
The functions $u_n|_{\Omega_n}$ define a holomorphic extension
of $u$ to $\C\setminus\mathcal A$, which we still denote by $u$.
The identity theorem also shows that the extended function satisfies
\eqref{eq:cohomology} on $\C\setminus(\mathcal A\cup P^{-1}(\mathcal A))$.
Equation \eqref{eq:iterated-cohomology} likewise holds on
$\C\setminus(\mathcal A\cup P^{-n}(\mathcal A))$
for every $n\ge1$.

\smallskip
\noindent\emph{Step 3: Removing the attracting periodic points.}
Fix $a\in\mathcal A$ of period $r$, and set $T:=P^r$ and
$\kappa:=\Lambda^r$. Then $T(a)=a$, $|T'(a)|<1<|\kappa|$, and
$
\kappa u(T(w))=T'(w)u(w)+R_r(w)
$
on a sufficiently small punctured neighborhood of $a$.
Choose $r_0>0$ and $\eta,c\in(0,1)$ such that the closed disk
$\overline{D(a,r_0)}$ contains no other point of $\mathcal A$,
$T^{-1}(a)\cap\overline{D(a,r_0)}=\{a\}$, and
$$
 |T(w)-a|\le\eta|w-a|,\qquad
 \frac{|T'(w)|}{|\kappa|}\le c,
 \qquad w\in\overline{D(a,r_0)}.
$$
Choose $0<\delta<r_0$ smaller than
$\min_{|w-a|=r_0}|T(w)-a|$. Every $v\in D(a,\delta)\setminus\{a\}$
has a preimage $w\in D(a,r_0)\setminus\{a\}$ under $T$, by the
argument principle. Such a preimage satisfies
$|w-a|\ge\eta^{-1}|v-a|$ and
$$
 |u(v)|\le c|u(w)|+M,\qquad
 M:=|\kappa|^{-1}\max_{\overline{D(a,r_0)}}|R_r|.
$$
Starting at $v$, choose successive preimages while their distance
from $a$ is less than $\delta$. The distances of these preimages from
$a$ increase by a factor at least $\eta^{-1}$, so after finitely many steps
the preimages reach the compact annulus $\delta\le|w-a|\le r_0$. If $B$ is
the maximum of $|u|$ on this annulus, iterating
$|u(v)|\le c|u(w)|+M$ along these preimages gives
$|u(v)|\le B+M/(1-c)$. Hence $u$ is bounded
near $a$ and the singularity of $u$ at $a$ is removable. This argument also
applies when $T'(a)=0$. Removing all points of $\mathcal A$
gives an entire extension, and \eqref{eq:cohomology} now holds
on $\C$.

\smallskip
\noindent\emph{Step 4: Affineness of the extension.}
Fix $1<\alpha<2$, and set
$N(r):=\max_{|v|\le r}|u(v)|/(1+|v|)^\alpha$.
There is a constant $C>0$ such that, for large $|v|$, every root
of $P(w)=v$ satisfies $|w|\le C|v|^{1/D}$. Applying \eqref{eq:cohomology} at such a root
gives
$$
 \frac{|u(v)|}{(1+|v|)^\alpha}
 \le C|v|^{-\frac{(D-1)(\alpha-1)}D}
 N(C|v|^{1/D})
 +C|v|^{\frac{D-2}D-\alpha}.
$$
Since $\alpha>1$ and $D\ge2$, both
$C|v|^{-\frac{(D-1)(\alpha-1)}D}$ and
$C|v|^{\frac{D-2}D-\alpha}$ tend to zero as $|v|\to\infty$. Choose $r_0$ so large
that, whenever $|v|\ge r_0$, the coefficient of $N$ is at most
$1/2$, the last term is at most $1$, and $C|v|^{1/D}\le|v|$.
For $r\ge r_0$, choose $v_r$ with $|v_r|\le r$ such that
$N(r)=|u(v_r)|/(1+|v_r|)^\alpha$.
If $|v_r|\ge r_0$, then
$N(r)\le\tfrac12N(C|v_r|^{1/D})+1\le\tfrac12N(r)+1$;
otherwise $N(r)\le N(r_0)$.
Thus $N(r)\le\max\{N(r_0),2\}$ for every $r\ge r_0$.
It follows that $u(w)=O((1+|w|)^\alpha)$.
For each $z\in\C$, Cauchy's estimate on circles centered at $z$
gives $|u''(z)|=O_z(r^{\alpha-2})$ as $r\to\infty$.
Since $\alpha<2$, we obtain $u''=0$, so $u(w)=Aw+B$.

\smallskip
\noindent\emph{Step 5: Coefficient comparison.}
The coefficients of $w^D$ and $w^{D-1}$ in
\eqref{eq:cohomology} give $(\Lambda-D)A=0$ and $Da_DB=0$,
respectively. Here the absence of a $w^{D-1}$ term in $P$ and the
degree bound on $R$ are both used. Hence $B=0$ and
$R=A(\Lambda P-wP')=A(DP-wP')$, as desired.
\end{proof}

\begin{remark}
The condition $\Lambda\ne D$ in the last assertion of
Lemma \ref{lem:extension} is necessary. For a hyperbolic polynomial
$P(w):=w^D+c$, the functions $u(w):=w$ and $R(w):=Dc$ satisfy
\eqref{eq:cohomology} with $\Lambda=D$.
\end{remark}

\section{From measurable history independence to a polynomial product}
In this section, we prove Theorem \ref{thm:rigidity} by extending the
measurable fast direction to a holomorphic line field near $K$ and
applying Lemma \ref{lem:extension} on periodic fibers.

Under \textup{(E$_0$)}, $DF$ is non-singular on $K$. The triangular
derivative cocycle $DF$ has exponents $\lambda_b>\lambda_v>0$, and its slow
Oseledec line is vertical. We write $\|\cdot\|$ for Euclidean vector norms and the induced
operator norms, $|\cdot|$ for complex modulus, and $\nabla$ for the
gradient in real coordinates. We use the description of $K$ and the
fiber Julia sets from \cite{Jonsson}.
Lemma \ref{lem:history-box} gives a local version of the history box
construction in \cite[Section 2.1]{BH}, under uniform expansion on $K$.

For a ball $U$ centered at $x\in K$, label the $d^2$ inverse
choices at each step by $1,\ldots,d^2$.
Write $\Sigma_U^-:=\{1,\ldots,d^2\}^{\N}$ for the space of
these sequences, and let $\beta_U$ be the product probability
measure assigning weight $d^{-2}$ to each choice.
For $t\in\Sigma_U^-$, write
$(x,x_{-1}(t),x_{-2}(t),\ldots)$ for the corresponding inverse
history over $x$.

\begin{lemma}
\label{lem:history-box}
Suppose that $F$ is uniformly expanding on $K$.
There exist $r>0$, $C_0>0$ and $\theta\in(0,1)$ such that,
for every $x\in K$ and $U:=B(x,r)$, each $t\in\Sigma_U^-$
determines holomorphic inverse branches $g_{t,n}$ of $F^n$
on $U$ satisfying $g_{t,n}(x)=x_{-n}(t)$.
These branches are compatible:
$F\circ g_{t,n+1}=g_{t,n}$ for $n\ge0$, where
$g_{t,0}:=\mathrm{id}_U$.
For $n\ge1$ and $y,y'\in U$, we have
\begin{equation}\label{eq:inverse-contraction}
 \|g_{t,n}(y)-g_{t,n}(y')\|
 \le C_0\theta^n\|y-y'\|.
\end{equation}
Moreover, the map $\Theta_U(y,t):=(y,g_{t,1}(y),g_{t,2}(y),\ldots)$
identifies $(K\cap U)\times\Sigma_U^-$ with $\pi_0^{-1}(U)$ and
satisfies
\begin{equation}\label{eq:local-product}
 \hat\mu|_{\pi_0^{-1}(U)}
 =(\Theta_U)_*\bigl((\mu|_U)\otimes\beta_U\bigr).
\end{equation}
\end{lemma}

\begin{proof}
We first construct the inverse branches. Choose $N$ so that $\|(DF_x^N)^{-1}\|<1/4$ on $K$.
Compactness and complete invariance of $K$ give a common $r>0$ such
that, for every $x\in K$ and $x_{-N}\in F^{-N}(x)$, the corresponding
inverse branch of $F^N$ satisfies
$$
 g\colon B(x,r)\to B(x_{-N},r),\qquad
 g(x)=x_{-N},\qquad \sup_{y\in B(x,r)}\|Dg(y)\|\le\tfrac12.
$$
Composing these inverse branches of $F^N$ gives inverse branches of
$F^{kN}$. After decreasing $r$, every inverse branch $h$ of $F^j$
on $B(x,r)$, with $x\in K$ and $0\le j<N$, satisfies
$\sup_{y\in B(x,r)}\|Dh(y)\|\le M$, with $M$ independent of $x,j,h$.
For $n=kN+j$, compose an inverse branch of $F^{kN}$ with
an inverse branch of $F^j$ based at the corresponding preimage.
The derivative of this composition is bounded by $M2^{-k}$.
Thus \eqref{eq:inverse-contraction} holds with
$\theta:=2^{-1/N}$ and $C_0:=M\theta^{-(N-1)}$.
These branches exhaust all inverse images over $K$, since $F^{-1}K=K$
and $DF$ is non-singular on $K$. Two inverse branches of $F^n$
that agree at one point of $U$ agree throughout $U$, by local uniqueness
and the identity theorem. Thus distinct branches remain distinct on $U$,
and the $d^2$ inverse choices at each step give the coding by $\Sigma_U^-$.
For $t\in\Sigma_U^-$, the branch $g_{t,n}$ is determined by
the first $n$ entries of $t$.

We now prove the local product formula. For a prefix $a$ of length $n$, let $g_a$ be the corresponding inverse
branch of $F^n$, and let $[a]\subset\Sigma_U^-$ consist of the sequences
beginning with $a$. For every Borel set $B\subset K\cap U$, the measure
of maximal entropy satisfies $\mu(g_a(B))=d^{-2n}\mu(B)$;
see \cite[Section 1.7]{DS}. Invariance of $\hat\mu$ therefore gives
$$
 \hat\mu\bigl(\Theta_U(B\times[a])\bigr)
 =\mu(g_aB)=d^{-2n}\mu(B)=\mu(B)\beta_U([a]).
$$
The sets $\Theta_U(B\times[a])$ generate the Borel $\sigma$-algebra
of $\pi_0^{-1}(U)$, proving \eqref{eq:local-product}.
\end{proof}

For fixed $t$, the ambient graph
$\Theta_U(U\times\{t\})$ is a full unstable plaque:
inverse iterates of any two points on this graph converge
exponentially.
By \eqref{eq:local-product}, the normalized conditional measures
obtained by fixing $t$ project to
$\mu_U:=\mu|_U/\mu(U)$ for $\beta_U$-almost every $t$.

\begin{proposition}\label{prop:ambient}
Under \textup{(E$_0$)}, suppose that the fast direction is
history-independent. Then there is a function $s$, holomorphic on a
neighborhood of $K$ and unique up to restriction to a smaller neighborhood,
such that $E^f(\hat x)=\Span(1,s(x_0))$ for $\hat\mu$-almost every $\hat x$.
On a sufficiently small neighborhood of $K$, the function $s$ satisfies
\begin{equation}\label{eq:slope}
 p'(z)s(F(z,w))=q_z(z,w)+q_w(z,w)s(z,w).
\end{equation}
\end{proposition}

\begin{proof}
\emph{Step 1: Typical histories and bounded distortion.}
Set $\alpha:=q_z/p'$ and $\beta:=q_w/p'$.
The functions $\alpha$ and $\beta$ are holomorphic on a neighborhood of $K$, and $\beta$ is
bounded away from zero there.
For $\hat\mu$-almost every $\hat x$, the fast line is non-vertical,
so we define its slope $\sigma(\hat x)\in\C$ by
$E^f(\hat x)=\Span(1,\sigma(\hat x))$.
$DF$-invariance of $E^f$ gives
$
\sigma(\hat F\hat x)=\alpha(x_0)+\beta(x_0)\sigma(\hat x).
$
Set $\Delta:=\lambda_b-\lambda_v>0$. The Birkhoff ergodic theorem
applied to $\hat F^{-1}$ and $\log|\beta|\circ\pi_0$ gives, for
$\hat\mu$-almost every history,
\begin{equation}\label{eq:typical-product}
 \lim_{n\to\infty}\frac1n
 \log\left|\prod_{j=1}^n\beta(x_{-j})\right|=-\Delta.
\end{equation}
Let $\vartheta(\hat x)$ be the angle between $E^f(\hat x)$ and
the vertical line. The angle property of the Oseledec splitting gives
$n^{-1}\log\sin\vartheta(\hat F^{-n}\hat x)\to0$ for
$\hat\mu$-almost every $\hat x$.
Since $\sin\vartheta(\hat x)=(1+|\sigma(\hat x)|^2)^{-1/2}$,
we obtain $\log^+|\sigma(\hat F^{-n}\hat x)|=o(n)$.

Fix $x\in K$ and let $U=B(x,r)$, with $r$ as in
Lemma \ref{lem:history-box}. Choose a neighborhood $W$ of $K$ with
$\sup_{\xi\in W}\|\nabla\log|\beta(\xi)|\|\le M$.
Decrease $r$ so that all inverse images of $U$ and the straight segments
between $g_{t,j}(y)$ and $g_{t,j}(y_*)$, for $y,y_*\in U$, lie in $W$.
By \eqref{eq:local-product} and Fubini's theorem, there is
$T_U\subset\Sigma_U^-$ with $\beta_U(T_U)=1$ such that, for every
$t\in T_U$ and $\mu_U$-almost every $y$, the history
$\hat x=\Theta_U(y,t)$ satisfies \eqref{eq:typical-product},
$\log^+|\sigma(\hat F^{-n}\hat x)|=o(n)$, and
$E^f(\hat x)=\ell(y)$.
Fix $t\in T_U$ and one such point $y_*$.
The bound $\|\nabla\log|\beta|\|\le M$ on $W$, together with
\eqref{eq:inverse-contraction}, gives a constant $C_1$ such that
$$
 \left|\sum_{j=1}^n
  \log\frac{|\beta(g_{t,j}(y))|}{|\beta(g_{t,j}(y_*))|}\right|
 \le C_1\sum_{j=1}^n\theta^j\|y-y_*\|,
 \qquad y\in U.
$$

Since $\sum_{j=1}^n\theta^j\le\theta/(1-\theta)$,
equation \eqref{eq:typical-product} at $\Theta_U(y_*,t)$
implies that, for every $\varepsilon\in(0,\Delta)$,
there is $C_{t,\varepsilon}<\infty$ such that
$$
 \sup_{y\in U}\left|\prod_{j=1}^n\beta(g_{t,j}(y))\right|
 \le C_{t,\varepsilon} e^{-(\Delta-\varepsilon)n},
 \qquad n\ge1.
$$

\smallskip
\noindent\emph{Step 2: A holomorphic representative.}
For each fixed $t\in T_U$, the products
$\prod_{k=1}^{j-1}\beta(g_{t,k}(y))$ decay exponentially
as $j\to\infty$, uniformly for $y\in U$.
Since $\alpha$ is bounded on $W$, the series
\begin{equation}
\label{eq:slope-series}
s_t(y):=\sum_{j=1}^{\infty}\alpha(g_{t,j}(y))
\prod_{k=1}^{j-1}\beta(g_{t,k}(y))
\end{equation}
converges absolutely and uniformly on $U$, defining a holomorphic function $s_t$.
For $\mu_U$-almost every $y\in K\cap U$, iterating the identity
$\sigma(\hat F\hat x)=\alpha(x_0)+\beta(x_0)\sigma(\hat x)$
along the inverse history $\Theta_U(y,t)$ gives
$$
\begin{aligned}
 \sigma(\Theta_U(y,t))
 &=\sum_{j=1}^{n}\alpha(g_{t,j}(y))
              \prod_{k=1}^{j-1}\beta(g_{t,k}(y))+r_{t,n}(y),\\
 r_{t,n}(y)&:=\left(\prod_{j=1}^{n}\beta(g_{t,j}(y))\right)
 \sigma(\hat F^{-n}\Theta_U(y,t)).
\end{aligned}
$$
Equation \eqref{eq:typical-product} and
$\log^+|\sigma(\hat F^{-n}\Theta_U(y,t))|=o(n)$
imply $r_{t,n}(y)\to0$ for $\mu_U$-almost every $y$.
Hence $s_t(y)=\sigma(\Theta_U(y,t))$ for $\mu_U$-almost every $y$.
Since the fast direction is history-independent,
$\Span(1,s_t(y))=\ell(y)$ for $\mu_U$-almost every $y$.

\smallskip
\noindent\emph{Step 3: Uniqueness, gluing, and invariance.}
The measure $\mu$ assigns zero mass to every proper local analytic set;
see \cite{DS}. 
If $V$ is a connected open set meeting $K=\supp\mu$ and a holomorphic
function $h$ on $V$ vanishes $\mu$-almost everywhere on $K\cap V$,
then $h\equiv0$ on $V$.

We now glue the local representatives. Choose finitely many balls $U_i$ covering $K$ and a history
$t_i\in T_{U_i}$ for each ball, and write $s_i:=s_{t_i}$.
On every overlap $U_i\cap U_j$ meeting $K$, we have $s_i-s_j=0$
$\mu$-almost everywhere. Since $U_i\cap U_j$ is connected and meets $K$,
the holomorphic function $s_i-s_j$ vanishes identically on $U_i\cap U_j$.
Choose smaller balls $V_i$ with $\overline V_i\subset U_i$ that still
cover $K$. If $\overline V_i\cap\overline V_j\cap K\ne\emptyset$,
then $U_i\cap U_j$ is connected and meets $K$, so $s_i=s_j$ on $U_i\cap U_j$.
For pairs $(i,j)$ with $\overline V_i\cap\overline V_j\cap K=\emptyset$,
remove the closed sets $\overline V_i\cap\overline V_j$ from
$\bigcup_i V_i$. The functions $s_i$ then glue to a holomorphic function
$s$ on the resulting neighborhood of $K$.
Two representatives agree $\mu$-almost everywhere near $K$, so their
difference vanishes identically near $K$, proving uniqueness.

Finally, we prove \eqref{eq:slope}. Since $F(K)=K$, both $s$ and
$s\circ F$ are defined near $K$. On this neighborhood, define the holomorphic function
$$
 \Psi(z,w):=p'(z)s(F(z,w))-q_z(z,w)-q_w(z,w)s(z,w).
$$
$DF$-invariance of the measurable fast line field gives $\Psi=0$
$\mu$-almost everywhere on $K$. Since a holomorphic function that is not identically zero on a
connected open set has a zero set of $\mu$-measure zero, $\Psi$ vanishes
on each connected component meeting $K$. Restricting to these components proves \eqref{eq:slope}.
\end{proof}

We call a periodic point $z_0$ of $p$ \emph{non-resonant} if
$|(p^n)'(z_0)|>d^n$ for some $n\ge1$ with $p^n(z_0)=z_0$.

\begin{lemma}\label{lem:periodic}
Under \textup{(E$_0$)}, the base polynomial is expanding on $J_p$.
For every $z_0\in J_p$ with $p^n(z_0)=z_0$, the return polynomial
$P(w):=Q_n(z_0,w)$, where $F^n(z,w)=(p^n(z),Q_n(z,w))$,
is hyperbolic and satisfies $\{z_0\}\times J_P\subset K$.
Moreover, the non-resonant periodic points are dense in $J_p$
and form a $p$-invariant set.
\end{lemma}

\begin{proof}
\emph{Step 1: Expansion of the base and periodic fibers.}
The projection of $K$ onto the $z$-coordinate is $J_p$, and each fiber Julia set over $J_p$
is contained in $K$; see \cite{Jonsson}.
Write $DF_x^m=\left(\begin{smallmatrix}A&0\\B&C\end{smallmatrix}\right)$.
Since $DF_x^m(1,-B/C)=(A,0)$, uniform expansion gives
$|A|\ge c\rho^m\|(1,-B/C)\|\ge c\rho^m$.
Thus $p$ is expanding on $J_p$.
Over a periodic point $z_0$, the fiber Julia set is $J_P$.
For every $w\in J_P$ and $k\ge1$, we have
$DF^{nk}_{(z_0,w)}(0,1)=(0,(P^k)'(w))$.
Since $(z_0,w)\in K$, uniform expansion of $F$ gives
$|(P^k)'(w)|\ge c\rho^{nk}$.
Hence $P$ is hyperbolic.

\smallskip
\noindent\emph{Step 2: Non-resonant periodic points.}
The projection of $\mu$ to the base is $\mu_p$, the measure of
maximal entropy of $p$.
The polynomial skew product Lyapunov formula gives
$$
 \lambda_v=\log d+
 \int_{J_p}\sum_{q_w(z,c)=0}G_z(c)\,d\mu_p(z)\ge\log d,
$$
where $G_z$ is the vertical escape rate function and critical points
are counted with multiplicity; see \cite{Jonsson} and
\cite[equation (3)]{AB}. Therefore $\lambda_b>\log d$.

We now prove that the non-resonant periodic points are dense in $J_p$. Fix a non-empty relatively open set $O\subset J_p$.
Choose a non-periodic recurrent point $z\in O$ such that
$\frac1n\sum_{k=0}^{n-1}\log|p'(p^k(z))|\to\lambda_b$.
Such points have full $\mu_p$-measure, and $\mu_p(O)>0$.
Choose $n_j\to\infty$ with $p^{n_j}(z)\to z$.
Choose a neighborhood $W$ of $J_p$ such that
$\sup_{\zeta\in W}\|\nabla\log|p'(\zeta)|\|\le M$.
Uniform expansion of $p$ on $J_p$ gives $r,C>0$ and $0<\theta<1$
such that the inverse branch $g_j$ of $p^{n_j}$ along the orbit of $z$
is defined on $D(p^{n_j}(z),r)$, with $g_j(p^{n_j}(z))=z$ and
$\operatorname{Lip}(g_j)\le C\theta^{n_j}$.
Take $r$ sufficiently small that, for every $j$ and
$0\le k\le n_j$, the set
$p^k(g_j(D(p^{n_j}(z),r)))$ has convex hull contained in $W$.
For large $j$, we have $g_j(\overline{D(z,r/2)})\subset D(z,r/2)$.
By the contraction mapping theorem, $g_j$ has a unique fixed point
$z_j$ in $D(z,r/2)$. The point $z_j$ satisfies $p^{n_j}(z_j)=z_j$ and
$|(p^{n_j})'(z_j)|>1$; hence $z_j\in J_p$, and $z_j\to z$.
The inverse branches defining $g_j$ give
$|p^k(z_j)-p^k(z)|\le C'\theta^{n_j-k}$ for $0\le k\le n_j$,
where $C'>0$ is independent of $j,k$.
The gradient bound gives
$\bigl|\log|p'(p^kz_j)|-\log|p'(p^kz)|\bigr|\le MC'\theta^{n_j-k}$.
Summing over $0\le k<n_j$ gives
$$
 \log|(p^{n_j})'(z_j)|
 =\sum_{k=0}^{n_j-1}\log|p'(p^kz)|+O(1).
$$
After division by $n_j$, the right side tends to $\lambda_b>\log d$.
Thus $z_j\in O$ is non-resonant for large $j$.
Since $O$ was arbitrary, the non-resonant periodic points are dense in $J_p$,
and hence infinite.

Finally, we show that the set of non-resonant periodic points is $p$-invariant.
If $p^n(z_0)=z_0$, then $(p^n)'(p^k(z_0))=(p^n)'(z_0)$ for
$0\le k<n$, since the factors in the derivative product are cyclically
permuted.
\end{proof}

\begin{proof}[Proof of Theorem \ref{thm:rigidity}]
We first prove that history independence implies \eqref{eq:product-form}.

\smallskip
\noindent\emph{Step 1: Centering the fiber polynomials.}
Since $F$ extends holomorphically to $\PP^2$, we can write
$q(z,w)=a_dw^d+\sum_{j=0}^{d-1}a_j(z)w^j$, where $a_d\ne0$ is
constant and $\deg a_j\le d-j$.
Write $b(z):=a_{d-1}(z)$ and $h(z):=b(z)/(da_d)$.
Since $\deg b\le1$, write $h(z)=cz+e$. Conjugating by
$H(z,w):=(z,w+h(z))$ replaces $q$ by
$
\widetilde q(z,v):=q(z,v-h(z))+h(p(z)),
$
whose $v^{d-1}$ coefficient vanishes, so the fiber polynomial is centered.
For $\widetilde F:=H\circ F\circ H^{-1}$,
we have $D\widetilde F^m_{H(x)}=DH\,DF_x^m(DH)^{-1}$.
Since $DH$ is constant and invertible, $\widetilde F$ is uniformly expanding
on $H(K)$, with unchanged Lyapunov exponents for $H_*\mu$.
Since $h(z)=cz+e$, the derivative $DH$ sends $(1,\sigma)$ to
$(1,\sigma+c)$, preserving history independence.
It therefore suffices to prove that the centered skew product is a
polynomial product. We now assume
$$
 q(z,w)=a_dw^d+\sum_{j=0}^{d-2}a_j(z)w^j.
$$
Proposition \ref{prop:ambient} gives a function $s$, holomorphic on a
neighborhood of $K$, satisfying \eqref{eq:slope}.

\smallskip
\noindent\emph{Step 2: Vanishing on non-resonant periodic fibers.}
Fix a non-resonant periodic point $z_0$ of period $n$, and set
$P(w):=Q_n(z_0,w)$, $D:=d^n$, and $\Lambda:=(p^n)'(z_0)$.
Composing the centered fiber polynomials shows that $Q_n(z,w)$ has
constant leading coefficient $a_d^{1+d+\cdots+d^{n-1}}$ and no $w^{D-1}$ term.
Thus $R(w):=\partial_zQ_n(z_0,w)$ has degree at most $D-2$.

Since $F^j(\{z_0\}\times J_P)\subset K$ for $0\le j<n$,
iterating \eqref{eq:slope} gives
$\Lambda s(z_0,P(w))=P'(w)s(z_0,w)+R(w)$ for $w$ near $J_P$.
Lemma \ref{lem:periodic} gives hyperbolicity of $P$, and
the choice of $z_0$ gives $|\Lambda|>D$.
Therefore by Lemma \ref{lem:extension}, we have $s(z_0,\cdot)=0$
near $J_P$.

\smallskip
\noindent\emph{Step 3: The product form.}
Every point of the base periodic orbit is non-resonant by
Lemma \ref{lem:periodic}. Thus $s$ vanishes
near the fiber Julia sets over both $z_0$ and $p(z_0)$.
Since $q(z_0,\cdot)$ maps $J_P$ onto the fiber Julia set over $p(z_0)$,
equation \eqref{eq:slope} gives $q_z(z_0,w)=0$ near $J_P$, and hence
for every $w\in\C$. There are infinitely many non-resonant $z_0$ by
Lemma \ref{lem:periodic}. Each coefficient in $w$ of the polynomial
$q_z(z,w)$ therefore vanishes at infinitely many $z_0$ and is zero.
Therefore the centered map is a polynomial product.

Returning to the original coordinates, let $Q_*$ be the centered fiber polynomial. Since
$H\circ F\circ H^{-1}(z,v)=(p(z),Q_*(v))$, we have
$q(z,w)=Q_*(w+cz+e)-cp(z)-e$. Setting $a:=-c$ and
$Q(v):=Q_*(v+e)-e$ gives \eqref{eq:product-form}.

\smallskip

We now prove that \eqref{eq:product-form} implies history independence.
The shear map $H_a$ conjugates $F$ to $(p,Q)$. Since $\lambda_b>\lambda_v$, the
fast Oseledec line of this product is horizontal for
$\mu_p\otimes\mu_Q$-almost every point, where $\mu_Q$ is the
measure of maximal entropy of $Q$.
Pulling the horizontal line field back by $H_a$ gives
$E^f=\Span(1,a)$ $\hat\mu$-almost everywhere,
proving history independence.
\end{proof}

\begin{remark}
The passage from \textup{(E)} to \textup{(E$_0$)} enlarges the
class of maps. For example, $(z,w)\mapsto(z^2-3,w^2-5/2)$
satisfies \textup{(E$_0$)} but not \textup{(E)}.
We first verify \textup{(E$_0$)}. The critical point $0$ of each
polynomial $z^2-3$ and $w^2-5/2$ escapes to infinity, so their product
is uniformly expanding on $K$. For $P_c(z):=z^2-c$, set $x_n:=P_c^n(0)$ and write $G_{P_c}$
for the escape rate function of $P_c$.
For $c=3,5/2$ and $n=4$, the escape rate series gives
$$
 2^{-n}\left(\log|x_n|-\frac{2c}{|x_n|^2}\right)
 \le G_{P_c}(0)\le2^{-n}\log|x_n|.
$$
Indeed, the subsequent $|x_j|$ increase, and
$-2t\le\log(1-t)\le0$ for $0\le t\le1/2$ bounds the
series terms. Since the Lyapunov exponent of $P_c$ equals $\log2+G_{P_c}(0)$ by \cite[equation (3)]{AB},
we obtain $1.1300<\lambda_b<1.1301$ and
$0.9979<\lambda_v<0.9980$.

We now show that \textup{(E)} fails.
At the fixed point $((1-\sqrt{13})/2,(1+\sqrt{11})/2)$,
the vertical multiplier has modulus $1+\sqrt{11}$, larger than the
modulus $\sqrt{13}-1$ of the base multiplier. A dominated splitting with the
vertical line slow is therefore impossible.
\end{remark}

\section{Binder--DeMarco equality in the base-dominant case}
\label{sec:base-dominant}

In this section, we prove Theorem \ref{thm:BD}.
We assume \textup{(E)} throughout and write
$\pi(z,w):=z$ and $L:=\log d$.

\begin{lemma}
\label{lem:local-geometry}
There is $r>0$ such that, for every $\hat x\in\hat K$, a
holomorphic flow box containing $B(x_0,2r)$ has coordinates $(z,v)$,
where $z$ is the original base coordinate and the fast plaques for that
history are the level sets of $v$.
Its plaques are graphs over the $z$-coordinate, with slopes uniformly
bounded on $B(x_0,r)$, and vary continuously with $\hat x$.
\end{lemma}

\begin{proof}
We use the history boxes $\Theta_U$ and inverse branches $g_{t,n}$ of
Lemma \ref{lem:history-box}, with history labels $t\in\Sigma_U^-$.
Set $\alpha:=q_z/p'$ and $\beta:=q_w/p'$.
The continuous splitting in \textup{(E)} on the compact set $\hat K$
makes the fast lines $E^f(\hat x)$ uniformly transverse to
$E^v(x_0)=\{0\}\times\C$.
Thus $\|DF^n|_{E^f(\hat x)}\|$ is uniformly comparable to
$|(p^n)'(\pi(x_0))|$, and \eqref{eq:domination} gives
$|\prod_{j=0}^{n-1}\beta(F^jx)|\le C_1\tau^n$ on $K$.
Choose $N \ge 1$ and a bounded neighborhood $W$ of $K$ such that
$\alpha$ and $\beta\circ F^j$, for $0\le j<N$, are holomorphic near
$\overline W$ and
$\sup_W|\prod_{j=0}^{N-1}\beta\circ F^j|<1/2$.
Shrink the balls $U$ so that all their inverse branch images lie in $W$.
For each $\ell\ge0$, the compatibility $F\circ g_{t,j+1}=g_{t,j}$ gives
$$
 \prod_{j=\ell N+1}^{(\ell+1)N}\beta(g_{t,j}(y))
 =\prod_{k=0}^{N-1}\beta\bigl(F^k(g_{t,(\ell+1)N}(y))\bigr).
$$
Since $g_{t,(\ell+1)N}(y)\in W$, the product on the right hand side
has modulus less than $1/2$. Writing
$n=mN+s$, with $0\le s<N$, and bounding the remaining $s$ factors gives
$|\prod_{j=1}^n\beta(g_{t,j}(y))|\le C_2 2^{-\lfloor n/N\rfloor}$, uniformly for $y\in U$, $t\in\Sigma_U^-$, and $n\ge1$.
Thus \eqref{eq:slope-series} defines uniformly bounded holomorphic
slopes $s_t$, continuous in the history. The strong unstable graph
transform identifies $\Span(1,s_t)$ with the fast direction, as in
\cite[proof of Proposition 2.1]{BH}.
Uniqueness of local solutions of $dw/dz=s_t(z,w)$ through a prescribed
point identifies the fast plaques with the integral curves of this
differential equation.
Taking $v$ to be the initial $w$-value on the vertical fiber through
$x_0$ gives the coordinates $(z,v)$.
Uniform ODE estimates and Cauchy estimates give the asserted flow boxes
and uniform slope bounds.
\end{proof}

By Lemma \ref{lem:history-box}, the normalized conditional measures on
fixed-history full unstable plaques project to $\mu_U:=\mu|_U/\mu(U)$,
as in \eqref{eq:local-product}.
Following \cite{QX} and \cite[Section 2.2]{BH}, let
$h_1:=H_{\hat\mu}(\xi^f\mid\hat F\xi^f)$, where $\xi^f$ is an
increasing measurable generator subordinate to the fast leaves.

\begin{proposition}\label{prop:dimensions}
The measure $\mu$ is exact-dimensional, with
\begin{equation}\label{eq:base-dimension}
 \Dim\mu=\frac{h_1}{\lambda_b}+\frac{2L-h_1}{\lambda_v}.
\end{equation}
In particular, $\Dim\mu=\BD(F)$ if and only if $h_1=L$.
\end{proposition}

\begin{proof}

Since $DF$ is non-singular on $K$,
$\log|\det D_{\mathbb R}F|\in L^1(\mu)$.
Since both Lyapunov exponents are positive, the full unstable
entropy is $h_2=h_\mu(F)=2L$.
Proposition 2.4 and equation (19) of \cite{QX}
give the dimension in \eqref{eq:base-dimension} for full unstable
conditionals. Lemma \ref{lem:history-box} identifies the projections of the
fixed-history full unstable plaque conditionals with $\mu_U$. On common refinements of the Qian--Xie full unstable partitions and the
fixed-history plaque partitions,
uniqueness of disintegration identifies the Qian--Xie full unstable
conditionals with the fixed-history plaque conditionals up to restriction
and normalization, as in \cite[proof of Proposition 2.3]{BH}.
Restrictions and normalizations on the common refinements preserve
almost-everywhere local dimensions. The projection $\pi_0$ is a local
diffeomorphism on each full unstable plaque, so it also preserves local
dimensions. Hence each $\mu_U$ is exact-dimensional with the dimension in
\eqref{eq:base-dimension}. A finite cover of $K$ by such balls gives
the same conclusion for $\mu$. Subtraction gives
$\Dim\mu-\BD(F)=(L-h_1)(1/\lambda_v-1/\lambda_b)$,
which proves the last assertion.
\end{proof}

Write $W^f_{\mathrm{loc}}(\hat x)$ for a local fast plaque projected
to $\C^2$. In the next two statements, these plaques are graphs over
$D(\pi(x_0),r_0)$ in the $z$-plane, with a common, sufficiently small
radius $r_0>0$.

\begin{lemma}
\label{lem:entropy-saturation}
We have $h_1\le L$. If $h_1=L$, then, for every $\hat x\in\hat K$, we have
\begin{equation}
\label{eq:fast-saturation}
 W^f_{\mathrm{loc}}(\hat x)\cap\pi^{-1}(J_p)\subset K.
\end{equation}
\end{lemma}

\begin{proof}
We adapt the entropy comparison argument of Buzzi--Fisher--Tahzibi
\cite[proof of Proposition 5.5]{BFT}. We compare the conditional probabilities
of $\hat\mu$ along fast leaves with measures obtained by transferring
$\mu_p$ to fast graphs through the $z$-projection.
Step 1 constructs an ambient partition whose atoms contain local fast disks
almost everywhere and whose restriction to $\hat K$ computes the fast entropy. Step 2 uses $h_1=L$ to show that the projected reference measures
are supported on $K$, yielding the inclusion on small disks for
$\hat\mu$-almost every history. Step 3 extends the inclusion to the
common local plaques for every history.

\smallskip
\noindent\emph{Step 1: An ambient subordinate partition.}

\smallskip
\noindent\emph{Step 1.1: Extending the fast graphs near $K$.}
We first choose an iterate whose inverse branches contract uniformly
near $K$. By uniform expansion and domination, choose an integer $N\ge1$
large enough that $G:=F^N=(P,Q_N)$ satisfies
$\|(DG)^{-1}\|<\kappa<1$ and $|\partial_wQ_N/P'|<1$ on a bounded
neighborhood $W_0$ of $K$. Here $P:=p^N$ and $D:=d^N$.
We now choose a smaller neighborhood $W$ of $K$ such that
\begin{equation}
\label{eq:ambient-neighborhood}
 G^{-1}(\overline W)\subset W,\qquad
 \sup_W\|(DG)^{-1}\|<\kappa,\qquad
 \sup_W|\partial_wQ_N/P'|<1.
\end{equation}
To obtain the first inclusion in \eqref{eq:ambient-neighborhood}, take
$W:=\{y:\operatorname{dist}(y,K)<\varepsilon\}$.
By properness of $G$ and $G^{-1}K=K$, choose $\varepsilon>0$ small enough
that $\overline W\subset W_0$ and $G^{-1}(\overline W)\subset W_0$.
For $y\in G^{-1}(\overline W)$, lift the segment from $Gy$ to a nearest
point of $K$ through $y$. The lift lies in $G^{-1}(\overline W)$ and ends
in $K$, so the inverse derivative bound gives
$\operatorname{dist}(y,K)\le\kappa\operatorname{dist}(Gy,K)
\le\kappa\varepsilon<\varepsilon$. Thus $y\in W$.

We next extend the fast graphs to histories whose current points may
lie outside $K$. Identify $\hat K$ with the natural extension of $G$ by retaining the
coordinates $x_{-jN}$. Let $\hat W$ be the space of $G$-inverse histories
in $W$, continue to write $\pi_0$ for the current coordinate projection
on $\hat W$, and let $T^{-1}$ delete the current coordinate. The forward
shift $T$ is defined when the forward image remains in $W$; on
$\hat K$ it is $\hat F^N$.
Apply Lemma \ref{lem:history-box} to $G$ to obtain history boxes on
coordinate balls, with $D^2$ inverse choices at each step. The inverse branches of $G^n$ contract by $\kappa^n$. Apply the slope series \eqref{eq:slope-series} to $G$, with
$\alpha:=\partial_zQ_N/P'$ and $\beta:=\partial_wQ_N/P'$.
The bound $\sup_W|\beta|<1$ in \eqref{eq:ambient-neighborhood}
gives uniform convergence,
so the construction in Lemma \ref{lem:local-geometry} defines invariant
holomorphic fast graphs on these ambient history boxes.

\smallskip
\noindent\emph{Step 1.2: Constructing the partition.}
To control how the partition boundaries cut fast graphs under inverse
iteration, choose a fine spatial grid partition $\mathcal P$ near $K$.
Write $\partial\mathcal P$ for the union of the cell boundaries. A translation of the grid can be chosen so that
\begin{equation}\label{eq:grid-summability}
 \sum_{j\ge1}\mu\{x:\operatorname{dist}(x,\partial\mathcal P)
 <\kappa^{j/2}\}<\infty.
\end{equation}
Indeed, the average of the $j$th term over translations
is $O(\kappa^{j/2})$. For each cell $R$, choose a coordinate ball $U$ with $\overline R\subset U$
on which the inverse branches and fast flow coordinates are defined, and
write $\Theta_U(y,t)$ for the history chart of $G$ over $U$.
Choose holomorphic transverse coordinates $v_t$, continuous in the
history label $t$, whose level sets are the fast plaques. Define a partition
$\mathcal Q$ by
$$
 \mathcal Q(\Theta_U(y,t))
 :=\{\Theta_U(y',t):y'\in R,\ v_t(y')=v_t(y)\},\qquad y\in R.
$$
Outside the union of the history boxes over these cells, set
$\mathcal Q(\hat x):=\{\hat x\}$.

We refine $\mathcal Q$ by requiring points to remain in the same cell and
fast plaque, with the same history label, at every past time. This defines the measurable partition
$
 \Xi(\hat x):=\{\hat y\in\mathcal Q(\hat x):
 T^{-j}\hat y\in\mathcal Q(T^{-j}\hat x)\text{ for every }j\ge0\}.
$

\smallskip
\noindent\emph{Step 1.3: Showing that the atoms contain fast disks.}
We now show that $\Xi(\hat x)$ contains a fast disk for
$\hat\mu$-almost every $\hat x$. Set
$d_j(\hat x):=\operatorname{dist}(\pi_0T^{-j}\hat x,\partial\mathcal P)$.
Since $T_*\hat\mu=\hat\mu$, we have
$\hat\mu\{d_j<\kappa^{j/2}\}
=\mu\{x:\operatorname{dist}(x,\partial\mathcal P)<\kappa^{j/2}\}$.
The Borel--Cantelli lemma and \eqref{eq:grid-summability} therefore give,
for $\hat\mu$-almost every $\hat x$,
\begin{equation}\label{eq:past-boundary-distance}
 d_j(\hat x)\ge\kappa^{j/2}\qquad\text{for all sufficiently large }j.
\end{equation}
Equation \eqref{eq:grid-summability} also gives $\mu(\partial\mathcal P)=0$;
invariance of $\hat\mu$ under $T$ then gives $d_j(\hat x)>0$ for
every $j\ge0$, for $\hat\mu$-almost every $\hat x$.

Choose $r_*>0$ uniformly for $\hat x\in\hat K$ so that the fast graph
through $x_0$ is defined over $D(\pi(x_0),r_*)$.
Write $\hat x=\Theta_U(x_0,t)$ and let $w=\varphi(z)$ denote the fast
graph through $x_0$ for this history label. For $0<r\le r_*$, define
$$
 \widehat D^f_r(\hat x)
 :=\{\Theta_U((z,\varphi(z)),t):z\in D(\pi(x_0),r)\}.
$$
Thus the same label $t$ determines the past coordinates of every point
of the graph. Lemma \ref{lem:local-geometry} gives a uniform bound
$|\varphi'|\le M_0$.
The graph over a line segment from $\pi(x_0)$ to
$z\in D(\pi(x_0),r)$ is a curve in $\C^2$ of length at most
$C_*r$, where $C_*:=\sqrt{1+M_0^2}$.
Its image under $j$ inverse iterates has length at most $C_*\kappa^jr$.
Set
$$
 r(\hat x):=\min\left\{r_*,\frac{1}{2C_*}
 \inf_{j\ge0}\kappa^{-j}d_j(\hat x)\right\}.
$$
By \eqref{eq:past-boundary-distance} and the fact that $d_j(\hat x)>0$ for all $j$,
this radius $r(\hat x)$ is positive for $\hat\mu$-almost every $\hat x$ and is measurable.
At time $-j$, each inverse image of such a graph curve stays in the cell
of $\mathcal P$ containing $\pi_0T^{-j}\hat x$,
since $C_*\kappa^jr(\hat x)\le d_j(\hat x)/2$.
Uniqueness of inverse branches ensures that the inverse image of the
graph curve at time $-j$ has the same history label as $T^{-j}\hat x$.
Invariance of the fast graphs under these inverse branches ensures that
the transverse coordinate along the inverse image curve equals that of
$T^{-j}\hat x$.
It follows that
$T^{-j}\widehat D^f_{r(\hat x)}(\hat x)
 \subset\mathcal Q(T^{-j}\hat x)$ for $j\ge0$,
so $\widehat D^f_{r(\hat x)}(\hat x)\subset\Xi(\hat x)$.

\smallskip
\noindent\emph{Step 1.4: Computing the conditional entropy.}
Finally, set $\xi:=\Xi|_{\hat K}$, so that
$\xi(\hat x)=\Xi(\hat x)\cap\hat K$ for $\hat x\in\hat K$.
We verify that $\xi$ is an increasing generator subordinate to the
fast leaves, so its conditional entropy computes the fast entropy of $T$. Indeed,
$C_n(\hat x):=T^{-n}\Xi(T^n\hat x)\subset\Xi(\hat x)$,
and $\diam\pi_0 C_n(\hat x)\le C\kappa^n$, with $C$ independent
of $\hat x$ and $n$. The atom of $T^{-n}\xi$ containing $\hat x$ is
$C_n(\hat x)\cap\hat K$.
All $C_n(\hat x)$ lie in $\Xi(\hat x)$,
where $\pi_0$ is injective. Since $\diam(\pi_0 C_n(\hat x))\to0$, we have
$\bigcap_{n\ge0}(C_n(\hat x)\cap\hat K)=\{\hat x\}$. Since $\xi$ is subordinate to the fast leaves, its partial entropy
for $T=\hat F^N$ is $Nh_1$. The chain rule therefore gives
$H_{\hat\mu}(T^{-n}\xi\mid\xi)=nNh_1$.

\medskip
\noindent\emph{Step 2: Entropy and reference mass.}

\smallskip
\noindent\emph{Step 2.1: Defining the reference measures.}
We first define the measures which will be compared with the fast leaf
conditionals. On a local fast graph lifted to $\hat W$ using a fixed
history label, define the reference measure by
$m^f(B):=\mu_p((\pi\circ\pi_0)(B))$ for Borel subsets $B$ of that graph.
The map $\pi\circ\pi_0$ assigns to a history the $z$-coordinate of its
current point and is injective on the graph. Thus $m^f$ transfers
$\mu_p$ to the graph, and the reference measures agree on overlapping
portions of the same lifted fast leaf.
Set $a(\hat x):=m^f(\Xi(\hat x))$; this function is constant on
each atom of $\xi$.
The function $a$ is measurable by integration in the flow coordinates, and
$0<a(\hat x)\le1$ for $\hat\mu$-almost every $\hat x$, since $\Xi(\hat x)$ contains
a disk centered over $J_p=\supp\mu_p$.

\smallskip
\noindent\emph{Step 2.2: Computing the mean negative logarithm.}
We next compute the mean negative logarithm of the normalized reference
mass of $C_n(\hat x)$.
We use truncation to justify the mean identity without assuming that
$\log a$ is integrable. The inverse branch defining $C_n(\hat x)$ projects to a single inverse branch of $P^n$. The constant Jacobian $D$ of $\mu_p$ under $P$ gives
\begin{equation}
\label{eq:reference-ratio}
 q_n(\hat x):=\frac{m^f(C_n(\hat x))}{a(\hat x)}
 =D^{-n}\frac{a(T^n\hat x)}{a(\hat x)}.
\end{equation}
Here $g:=-\log q_1\ge0$, and telescoping gives
$n^{-1}\sum_{j=0}^{n-1}g(T^j\hat x)
=\log D+n^{-1}(\log a(\hat x)-\log a(T^n\hat x))$.
Since $T_*\hat\mu=\hat\mu$ and $\log a$ is finite $\hat\mu$-almost everywhere,
$n^{-1}(\log a(\hat x)-\log a(T^n\hat x))$ tends to zero in
$\hat\mu$-probability.
Choose integers $n_k\to\infty$ such that, for $\hat\mu$-almost every $\hat x$, we have
$$
 \frac1{n_k}\sum_{j=0}^{n_k-1}g(T^j\hat x)\to \log D.
$$
For each integer $M\ge1$, set $g_M:=\min(g,M)$.
By the Birkhoff ergodic theorem, there is a measurable function $g_M^*$ such
that, for $\hat\mu$-almost every $\hat x$, we have
$$
 \frac1n\sum_{j=0}^{n-1}g_M(T^j\hat x)\to  g_M^*(\hat x).
$$
Since $g_M\le g$, comparison along $n_k$ gives $0\le g_M^*\le\log D$
for $\hat\mu$-almost every $\hat x$. The truncated averages are bounded by $M$, so dominated
convergence and $T$-invariance give
$\int g_M\,d\hat\mu=\int g_M^*\,d\hat\mu\le\log D$.
Monotone convergence as $M\to\infty$ yields
$\int g\,d\hat\mu\le\log D<\infty$.
Conversely, Fatou's lemma along $n_k$ and $T$-invariance give
$$
 \log D\le\liminf_{k\to\infty}
 \int\frac1{n_k}\sum_{j=0}^{n_k-1}g(T^j\hat x)\,d\hat\mu
 =\int g\,d\hat\mu.
$$
Therefore $\int g\,d\hat\mu=\log D$.
Consequently, \eqref{eq:reference-ratio} gives
$\int-\log q_n\,d\hat\mu=n\log D$.

\smallskip
\noindent\emph{Step 2.3: Comparing entropy with reference mass.}
We now compare the reference masses with the conditional probabilities
on the fast leaves. Disintegrate $\hat\mu$ over $\xi$, writing the conditional probabilities
as $\hat\mu_A$ and the quotient measure as $\bar\mu$.
Fix $n\ge1$. For $\bar\mu$-almost every atom $A$, fix $\hat x\in A$.
Since $T$ preserves $\hat\mu$, for $\hat\mu_A$-almost every
$\hat y\in A$ the atom $\Xi(T^n\hat y)$ contains an ambient fast disk.
The pullback of this disk gives an ambient leaf neighborhood in
$C_n(\hat y)\subset\Xi(\hat x)$, so $m^f(C_n(\hat y))>0$.
Distinct sets $C_n(\hat y)$ are disjoint and have positive reference mass;
since $m^f(\Xi(\hat x))\le1$, at most countably many occur.
Denote the distinct sets $C_n(\hat y)$ with positive conditional probability by $C_{A,j}$, and set
$p_j:=\hat\mu_A(C_{A,j}\cap\hat K)>0$ and
$q_j:=m^f(C_{A,j})/a(\hat x)$.
For $\hat y\in C_{A,j}\cap\hat K$, we have $q_n(\hat y)=q_j$.
Then $\sum_jp_j=1$, whereas only $\sum_jq_j\le1$ is known,
since we have not shown that $\bigcup_jC_{A,j}$ has full reference
measure in $\Xi(\hat x)$.
The log sum inequality gives
\begin{equation}\label{eq:entropy-comparison}
 \sum_jp_j\log\frac{p_j}{q_j}
 \ge-\log\sum_jq_j\ge0.
\end{equation}
Integrating $\sum_jp_j\log(p_j/q_j)$ over the atoms $A$
with respect to $\bar\mu$ gives
$n\log D-H_{\hat\mu}(T^{-n}\xi\mid\xi)=nN(L-h_1)$.
Since $\sum_jp_j\log(p_j/q_j)\ge0$, we obtain $h_1\le L$.
If $h_1=L$, then $\sum_jp_j\log(p_j/q_j)=0$
for $\bar\mu$-almost every atom $A$.

\smallskip
\noindent\emph{Step 2.4: Deducing local saturation.}
We next show that $h_1=L$ forces the projected reference measure to be
supported on $K$. If $h_1=L$, equality in \eqref{eq:entropy-comparison} gives
$\sum_jq_j=1$ for $\bar\mu$-almost every atom $A$, for each $n\ge1$.
Intersecting these full-measure sets of atoms over all $n\ge1$ gives
a single set of full $\bar\mu$-measure on which the equality holds for every $n$.
Since $\sum_jq_j=1$, the union $\bigcup_jC_{A,j}$ has full reference
measure in $\Xi(\hat x)$. Each $C_{A,j}$ meets $\hat K$ and satisfies
$\diam(\pi_0 C_{A,j})\le C\kappa^n$, so every point of
$\pi_0 C_{A,j}$ is within $C\kappa^n$ of $K$.
Taking the intersection, over all $n\ge1$, of the unions
$\bigcup_jC_{A,j}$ of full $m^f$-measure and using closedness of $K$
gives $m^f(\Xi(\hat x)\setminus\pi_0^{-1}(K))=0$.
Complete invariance under $G$ then gives
$m^f(\Xi(\hat x)\setminus\hat K)=0$ for $\hat\mu$-almost every $\hat x$.

We convert this measure statement into an inclusion of graph points
over $J_p$. If a point of $\pi_0\widehat D^f_{r(\hat x)}(\hat x)$ over $J_p$
were outside $K$, closedness of $K$ would give a leaf neighborhood disjoint
from $K$ with positive reference mass, since $\supp\mu_p=J_p$.
Thus, for $\hat\mu$-almost every $\hat x$, we have
\begin{equation}\label{eq:small-fast-inclusion}
 \pi_0\widehat D^f_{r(\hat x)}(\hat x)\cap\pi^{-1}(J_p)\subset K.
\end{equation}

\medskip
\noindent\emph{Step 3: All local plaques.}
We first extend \eqref{eq:small-fast-inclusion} to plaques of a common
$z$-radius for $\hat\mu$-almost every history.
We work on a $T$-invariant set of full $\hat\mu$-measure where
\eqref{eq:small-fast-inclusion} holds at every iterate. The sets $\{r(\hat x)\ge1/m\}$ cover
this set up to a set of $\hat\mu$-measure zero.
By the Poincar\'e recurrence theorem applied to $T^{-1}$
on the sets $\{r\ge1/m\}$, for $\hat\mu$-almost every $\hat x$
there are an integer $m=m(\hat x)\ge1$ and integers $n_j\to\infty$
such that $r(T^{-n_j}\hat x)\ge1/m$. Take the common radius $r_0$ in the $z$-plane
to be at most $r_*$.
Choose $n=n_j$ large enough that
$C_*\kappa^nr_*<1/m\le r(T^{-n}\hat x)$.
By invariance of the fast graphs and contraction of the inverse branches,
the inverse branch of $G^n$ determined by $\hat x$ maps
$W^f_{\mathrm{loc}}(\hat x)$ into
$\pi_0\widehat D^f_{r(T^{-n}\hat x)}(T^{-n}\hat x)$, whose points
over $J_p$ lie in $K$.
If $y\in W^f_{\mathrm{loc}}(\hat x)\cap\pi^{-1}(J_p)$, its image
under the inverse branch of $G^n$ determined by $\hat x$ still
projects into $J_p$, since $P^{-1}J_p=J_p$, and hence lies in $K$.
Applying $G^n$ and using $G(K)=K$ gives $y\in K$.

We now extend \eqref{eq:fast-saturation} to every history.
Equation \eqref{eq:local-product} gives $\supp\hat\mu=\hat K$, so any history
can be approximated by histories for which \eqref{eq:fast-saturation} holds.
After a uniform reduction of $r_0$, the approximating fast graphs are
defined on the $z$-disk of the limiting plaque.
Keeping $z\in J_p$ in this disk fixed, continuity of the fast graphs
and closedness of $K$ give \eqref{eq:fast-saturation} for every history.
\end{proof}

\begin{proposition}\label{prop:saturation-rigidity}
If \eqref{eq:fast-saturation} holds, then the fast direction is
history-independent.
\end{proposition}

\begin{proof}
Write $K_z:=\{w:(z,w)\in K\}$.
For a periodic base point $z_0\in J_p$ of period $n$, set
$P_0(w):=Q_n(z_0,w)$.
Uniform expansion on $K$ implies vertical expansion over $J_p$.
By \cite[Proposition 3.5]{Jonsson}, $K_{z_0}$ is the fiber
Julia set over $z_0$.
Since the fiber iterates at times $nk$ are $P_0^k$,
we have $K_{z_0}=J_{P_0}$.
We use local symmetry rigidity when some periodic fiber Julia set is
non-smooth; otherwise, we prove the product form directly.

\smallskip
\noindent\emph{Case 1: A non-smooth periodic fiber.}
Suppose there exists a periodic base point $z_0\in J_p$ such that
$J_{P_0}$ is neither a circle nor a segment.

\smallskip
\noindent\emph{Step 1: Constructing holonomies that preserve the fibers.}
We compare fast plaques by following them between nearby vertical fibers.
Fix $w_0\in J_{P_0}$ and a history box $\Theta_U$ about $(z_0,w_0)$,
as in Lemma \ref{lem:history-box}.
For a history label $t$, let $(z,H_t(z,w))$ be the point over $z$ on
the fast graph determined by $\Theta_U((z_0,w),t)$.
Lemma \ref{lem:local-geometry} defines these holonomies as solutions of
$$
 \partial_zH_t(z,w)=s_t(z,H_t(z,w)),\qquad H_t(z_0,w)=w,
$$
where $s_t$ is the holomorphic fast slope for the history label $t$. The maps $H_t$ are holomorphic on
$D_b\times D_w$, where $D_b\ni z_0$ and $D_w\ni w_0$ are disks
independent of the history $t$.

Shrink the disks $D_b$ and $D_w$ uniformly in $t$ so that the
$H_t$-graph over the line segment from $z_0$ to $z$ lies in the local
fast plaques based at $(z_0,w)$ and $(z,H_t(z,w))$. We then have
$$
 w\in J_{P_0}\quad\Longleftrightarrow\quad H_t(z,w)\in K_z,
 \qquad z\in D_b\cap J_p,\ w\in D_w.
$$
For the forward implication, apply \eqref{eq:fast-saturation} at
$(z_0,w)\in K$. For the reverse implication, $F^{-1}K=K$ ensures that
$\Theta_U((z,H_t(z,w)),t)\in\hat K$ whenever $(z,H_t(z,w))\in K$.
The fast graph for this history and the graph defined by $H_t$
both solve $dh/dz=s_t(z,h)$ and pass through $(z,H_t(z,w))$.
Uniqueness therefore identifies the two graphs.
Applying \eqref{eq:fast-saturation} at $(z,H_t(z,w))$ and returning to
$z_0$ gives $w\in K_{z_0}=J_{P_0}$.

\smallskip
\noindent\emph{Step 2: Obtaining local symmetries.}
We now compare two history labels $t,u$ by going from $z_0$ to $z$
along the $t$-plaques and returning along the $u$-plaques.
Choose a disk $D'_w$ about $w_0$ with $\overline{D'_w}\subset D_w$.
For $w_1,w_2\in D_w$, if $H_u(z,w_1)=H_u(z,w_2)$, uniqueness of solutions to
$dh/dz=s_u(z,h)$ makes the two solutions agree at $z_0$, so $w_1=w_2$.
Thus each $H_u(z,\cdot)$ is injective on $D_w$.
The history label space $\Sigma_U^-$ is compact, the holonomies depend
continuously on the labels, and $H_t(z_0,\cdot)=\mathrm{id}$.
After shrinking $D_b$ uniformly in $t,u$, we therefore have
$H_t(z,D'_w)\subset H_u(z,D_w)$ for every $z\in D_b$.
Thus $\sigma_z:=H_u(z,\cdot)^{-1}\circ H_t(z,\cdot)$ is well-defined
on $D'_w$, using the inverse of $H_u(z,\cdot)$ on $D_w$.
The maps $\sigma_z$ depend holomorphically on $z\in D_b$ and satisfy
$\sigma_{z_0}=\mathrm{id}$.
For $z\in D_b\cap J_p$ and $w\in D'_w$, the identity
$H_u(z,\sigma_z(w))=H_t(z,w)$ gives
$$
 w\in J_{P_0}\quad\Longleftrightarrow\quad H_t(z,w)\in K_z
 \quad\Longleftrightarrow\quad \sigma_z(w)\in J_{P_0}.
$$
Thus, for every $z\in D_b\cap J_p$, we have
$\sigma_z^{-1}(J_{P_0})\cap D'_w=J_{P_0}\cap D'_w$,
so $\sigma_z$ is a local symmetry of $J_{P_0}$.

\smallskip
\noindent\emph{Step 3: Proving local history independence.}
We use local symmetry rigidity to show that $\sigma_z=\mathrm{id}$. Levin's local symmetry theorem \cite{Levin}, in the form of
\cite[Theorem 2.1]{DFG}, implies that distinct local symmetries of a
non-smooth Julia set cannot converge to the identity on a fixed domain.
Since $J_p$ has no isolated points, take distinct $z_j\in D_b\cap J_p$
tending to $z_0$. Then $\sigma_{z_j}\to\mathrm{id}$ locally uniformly on $D'_w$.
Since $\sigma_{z_j}\to\mathrm{id}$, each fixed map different from
$\mathrm{id}$ can occur only finitely often in this sequence.
If $\sigma_{z_j}\ne\mathrm{id}$ for infinitely many $j$, we could therefore
extract a sequence of distinct local symmetries converging to the identity,
contradicting Levin's theorem.
Hence $\sigma_{z_j}=\mathrm{id}$ for all sufficiently large $j$.
For each fixed $w\in D'_w$, the identity theorem applied to
$z\mapsto H_t(z,w)-H_u(z,w)$ gives $H_t=H_u$ on $D_b\times D'_w$.
At $(z,H_t(z,w))$, the fast line is spanned by
$(1,\partial_zH_t(z,w))$. Since $t,u$ were arbitrary, write $H$ for the common holonomy and set
$V:=\{(z,H(z,w)):z\in D_b,\ w\in D'_w\}$.
All histories determine the same fast line on $V\cap K$.
The map $(z,w)\mapsto(z,H(z,w))$ has non-zero Jacobian, since each
$H(z,\cdot)$ is locally biholomorphic. Hence $V$ is open.
Since $(z_0,w_0)\in V\cap K$ and $K=\supp\mu$, we have $\mu(V)>0$.

\smallskip
\noindent\emph{Step 4: Extending history independence by ergodicity.}
We extend history independence from $V\cap K$ to $\mu$-almost every
current point. Disintegrate $\hat\mu=\int\hat\mu_x\,d\mu(x)$ over $\pi_0$ and set
$\eta_x:=(E^f)_*\hat\mu_x$. This is the distribution of fast directions
among histories with current point $x$; it is Dirac precisely when the
direction is constant for $\hat\mu_x$-almost every history.
The inverse branch weights give
$\eta_y=d^{-2}\sum_{Fx=y}(\mathbb P DF_x)_*\eta_x$
for $\mu$-almost every $y$, where $\mathbb P DF_x$ is the induced map
on complex lines; see \cite[Proposition 3.1]{BH}.
If $\eta_y$ is Dirac, each measure $(\mathbb P DF_x)_*\eta_x$ in this
positive-weight average is supported on the same single line as $\eta_y$. Since $DF_x$ is invertible,
each $\eta_x$ in the sum is therefore Dirac. Hence the set $A:=\{x:\eta_x\text{ is Dirac}\}$
satisfies $F^{-1}A\subset A$ modulo $\mu$.
Since $\mu(F^{-1}A)=\mu(A)$, we have $\mu(A\triangle F^{-1}A)=0$.
Ergodicity of $F$ with respect to $\mu$ gives $\mu(A)\in\{0,1\}$.
Since $\mu(A)\ge\mu(V)>0$, we obtain $\mu(A)=1$.

\smallskip
\noindent\emph{Case 2: All periodic fiber Julia sets are smooth.}
In this case we prove the product form directly.
The smooth Julia sets of polynomials are circles or segments;
see \cite[Section 2.2]{DFG}. Hyperbolicity of the return polynomials $P_0$, given by
Lemma \ref{lem:periodic}, excludes segments, so every periodic fiber
Julia set is a circle.
Apply the affine centering change from Step 1 of the proof of
Theorem \ref{thm:rigidity}, which preserves history independence, and
retain the notation $q$ for the transformed fiber polynomial.
Write $a_d\ne0$ for its constant leading coefficient.
Each fiber polynomial has zero coefficient of $w^{d-1}$; composing them
shows that every return polynomial is centered as well.
A polynomial with a circle Julia set is affinely conjugate to a
monomial. Since $Q_n(z_0,\cdot)$ is centered, the affine conjugacy
has zero translation term. Hence $Q_n(z_0,w)=A_nw^{d^n}$ at every
periodic $z_0\in J_p$,
where $A_n:=a_d^{1+d+\cdots+d^{n-1}}$.

We next deduce that each single step fiber polynomial over a periodic
point is a monomial. Write $q_{z_0}(w):=q(z_0,w)$.
If $n=1$, then $q_{z_0}(w)=Q_1(z_0,w)=a_dw^d$.
If $n\ge2$, set $R(v):=Q_{n-1}(p(z_0),v)$, so
$Q_n(z_0,\cdot)=R\circ q_{z_0}$.
Since $R$ is a composition of centered fiber polynomials, it is centered
of degree $d^{n-1}$.
The chain rule gives $(R\circ q_{z_0})'=(R'\circ q_{z_0})q_{z_0}'$.
Since $(R\circ q_{z_0})'$ vanishes only at zero, so does $q_{z_0}'$.
Thus $q_{z_0}(w)=a_dw^d+b$, and
$R(a_dw^d+b)=A_nw^{d^n}$.
Set $v=a_dw^d+b$ and use $w^{d^n}=(w^d)^{d^{n-1}}$.
Since $v$ ranges over $\C$, we obtain
$R(v)=c(v-b)^{d^{n-1}}$, where $c:=A_n/a_d^{d^{n-1}}\ne0$.
The coefficient of $v^{d^{n-1}-1}$ is $-cd^{n-1}b$; since $R$ is centered,
this coefficient is zero, so $b=0$.
Finally, each coefficient of $q(z,w)-a_dw^d$ as a polynomial in $w$
is a polynomial in $z$. Each coefficient vanishes at the infinitely many
periodic points in $J_p$, hence identically. Thus the centered map is
$(p(z),a_dw^d)$, whose fast direction is horizontal since
$\lambda_b>\lambda_v$. The affine centering change preserves history
independence, which proves the conclusion in the original coordinates.
\end{proof}

\begin{proof}[Proof of Theorem \ref{thm:BD}]
Theorem \ref{thm:rigidity} identifies history independence with
\eqref{eq:product-form}, already under \textup{(E$_0$)}.
If $F$ has the product form \eqref{eq:product-form}, the shear map $H_a$
carries $\mu$ to $\mu_p\otimes\mu_Q$.
The Ma\~n\'e--Manning formula \cite{Mane1988,Manning1984} gives exact
dimensions $L/\lambda_b$ and $L/\lambda_v$ for the factors.
Local dimensions add under products, and affine conjugacy preserves
dimension on compact supports, so $\Dim\mu=\BD(F)$.
Conversely, under \textup{(E)}, Proposition \ref{prop:dimensions}
turns dimension equality into $h_1=L$.
Lemma \ref{lem:entropy-saturation} and
Proposition \ref{prop:saturation-rigidity} then give history independence.
Exact dimensionality follows from Proposition \ref{prop:dimensions}.
\end{proof}

\begin{remark}
\label{rem:boundary}
The boundary $\lambda_v=\log d$ admits both equality $\Dim\mu_{F}=\BD(F)$ and strict
inequality $\Dim\mu_{F}>\BD(F)$. For example, let
$F_t(z,w):=(z^2-100,w^2+tz)$, where $|t|\le1/100$.
On $J_p$, $|z|\le11$ and $|p'(z)|\ge B:=2\sqrt{89}$.
For $z\in J_p$, the disks $|w|\le3/4$ are strictly forward invariant,
while points with $|w|\ge5/4$ escape. The fiber iterates are uniformly
bounded on $|w|<3/4$, so these open disks lie in the fiber Fatou sets. Hence $K$ lies in
$3/4\le|w|\le5/4$, and all vertical critical orbits over $J_p$ are bounded.
The Lyapunov formula of \cite{Jonsson} gives $\lambda_v=\log2$.
On $K$, the inverse derivative has operator norm, for the maximum
norm on $\C^2$, at most
$\max\{B^{-1},\frac23(1+|t|/B)\}<1$, and
$|q_w/p'|\le5/(2B)<1$. The slope graph transform therefore gives
\textup{(E)}. For $t\ne0$, comparison of the $w$ coefficient in
\eqref{eq:product-form} forces $a=0$, which is incompatible with $tz$.
Together with \eqref{eq:base-dimension} and the bound $h_1\le L$ in
Lemma \ref{lem:entropy-saturation},
Theorem \ref{thm:BD} gives $\Dim\mu_{F_t}>\BD(F_t)$ for $t\ne0$,
whereas $F_0$ is a product and satisfies equality.
Moreover, $F_0$ satisfies \textup{(E)} but not \textup{(H)} of \cite{BH}:
the support $J_{z^2-100}\times S^1$ of $\mu_{F_0}$ contains circles, whereas
\textup{(H)} gives a totally disconnected support.
\end{remark}

\appendix
\section{The fiber-dominant case}
\label{sec:fibered}
We explain how the Ledrappier--Young dimension argument in
\cite[proof of Proposition 5.3]{BH} gives the following result
without uniform expansion or domination assumptions.
The dimension formula in this generality was also proved
independently in \cite{Genes} using a fibered normal form theorem.

\begin{proposition}\label{prop:fiber-dominant}
Suppose that $\lambda_v>\lambda_b$.
Then the fast direction is vertical and hence history-independent.
Moreover, $\mu_F$ is exact-dimensional and
$\Dim\mu_F=\BD(F)$.
\end{proposition}

For the rest of this appendix, assume $\lambda_v>\lambda_b$. We first show that the  Ledrappier--Young formula of
\cite[Proposition 2.3]{BH} remains valid without hypothesis \textup{(H)}, and then
apply the dimension argument in \cite[proof of Proposition 5.3]{BH}.
The following lemma gives the required local product structure.

\begin{lemma}\label{lem:nonuniform-history-boxes}
Fix $0<a<\lambda_b$. Up to a set of $\hat\mu$-measure zero,
$\hat K$ is covered by countably many measurable boxes
$\mathcal B=\Theta((K\cap U)\times S)$, where $U$ is a coordinate
ball and
$
 \Theta(x,s)=(x,g_{s,1}(x),g_{s,2}(x),\ldots).
$
The compatible holomorphic inverse branches $g_{s,n}$ of $F^n$
are defined on $U$ and satisfy
$\operatorname{Lip}(g_{s,n}|_U)\le Ce^{-an}$, with $C$ depending
on the box. There is a finite measure $\beta$ on $S$ such that
$
 \hat\mu|_{\mathcal B}
 =\Theta_*\bigl((\mu|_U)\otimes\beta\bigr).
$
If $\hat\mu(\mathcal B)>0$, the conditional probabilities obtained
by fixing $s$ project to $\mu_U:=\mu|_U/\mu(U)$ for
$\beta$-almost every $s$.
\end{lemma}

\begin{proof}
For $\hat\mu$-almost every history, the inverse branches along that
history are defined on a common neighborhood of $x_0$, with
Lipschitz constants bounded by $C(\hat x)e^{-an}$; see
\cite[Theorem 1.4 (3)]{BDM}.
Restricting to a ball from a countable basis and increasing
$C(\hat x)$ to an integer gives countably many choices of $U,C$.

Fix $U,C$. We form a tree whose level-$n$ vertices are tuples
$(g_1,\ldots,g_n)$ of holomorphic inverse branches on $U$ satisfying
$F\circ g_j=g_{j-1}$, with $g_0=\mathrm{id}_U$, and
$\operatorname{Lip}(g_j|_U)\le Ce^{-aj}$ for $1\le j\le n$.
A child extends the tuple by one branch, so each vertex has at most
$d^2$ children. Let $S$ be the compact space of infinite paths.
Evaluation defines an injective Borel map $\Theta$ on
$(K\cap U)\times S$. Set $\mathcal B:=\Theta((K\cap U)\times S)$.
These boxes cover $\hat\mu$-almost every history by the inverse branch
theorem. For fixed $s$, the ambient graph $\Theta(U\times\{s\})$
is a local full unstable plaque, since its inverse iterates contract
exponentially.

We now construct $\beta$ and prove the product identity.
Since the Jacobian of $\mu$ is the constant $d^2$, for every inverse
branch $g_n$ of $F^n$ on $U$ and every Borel set $A\subset K\cap U$, we have
$
 \hat\mu\{\hat x:x_0\in A,\ x_{-n}=g_n(x_0)\}
 =\mu(g_n(A))=d^{-2n}\mu(A);
$
see \cite[Section 1.7]{DS}.
For a finite path $v$, let $[v]\subset S$ consist of the
infinite paths extending $v$, and let $N_m(v)$ count the finite
extensions of $v$ to level $m$. The union of the corresponding cylinders
over $A$ has measure $d^{-2m}N_m(v)\mu(A)$, and these unions
decrease to $\Theta(A\times[v])$ as $m\to\infty$. Thus
$
 \beta([v]):=\lim_{m\to\infty}d^{-2m}N_m(v)
$
defines consistent cylinder masses, since each vertex has finitely
many children. These cylinder masses determine a finite measure on $S$.
Continuity from above gives
$\hat\mu(\Theta(A\times[v]))=\mu(A)\beta([v])$.
The sets $\Theta(A\times[v])$ form a generating $\pi$-system
on $\mathcal B$, so uniqueness of finite measures gives
the product identity.

Finally, if $\hat\mu(\mathcal B)>0$, normalizing and disintegrating
the product identity with respect to $s$ gives projected conditionals
$\mu_U$ for $\beta$-almost every $s$.
\end{proof}

Write $\pi(z,w):=z$, so that $\pi_*\mu=\mu_p$, and disintegrate
$\mu=\int\mu_z\,d\mu_p(z)$ along the vertical fibers.
Set $L:=\log d$ and
$h_1:=H_{\hat\mu}(\xi^f\mid\hat F\xi^f)$, where $\xi^f$ is an
increasing generator subordinate to the fast leaves, as in
\cite[Section 2.2]{BH}. Here the fast direction is vertical.

\begin{proposition}\label{prop:fibered-dimensions}
The measure $\mu$ is exact-dimensional, and $\mu_z$ is
exact-dimensional for $\mu_p$-almost every $z$, with dimensions
$$
 \Dim\mu=\frac{h_1}{\lambda_v}+\frac{2L-h_1}{\lambda_b},
 \qquad
 \Dim\mu_z=\frac{h_1}{\lambda_v}.
$$
For every positive measure box in
Lemma \ref{lem:nonuniform-history-boxes}, the transverse measure
$\zeta_U:=\pi_*\mu_U$ is exact-dimensional of dimension
$(2L-h_1)/\lambda_b$.
\end{proposition}

\begin{proof}
By \cite{DS}, $\log|\det D_{\mathbb R}F|\in L^1(\mu)$.
Since both Lyapunov exponents are positive, the full unstable
entropy is $h_2=h_\mu(F)=2L$.
Proposition 2.4 and equation (19) of \cite{QX} give exact dimensions
$\delta_f=h_1/\lambda_v$ and
$\delta_u=h_1/\lambda_v+(2L-h_1)/\lambda_b$
for the conditional measures on the fast and full unstable
plaques, respectively.
Since the vertical direction is invariant and has the larger
exponent, the local strong unstable plaques are vertical disks.

We next prove the dimension formulas for $\mu$ and $\mu_z$.
Lemma \ref{lem:nonuniform-history-boxes} identifies the projections
of the fixed-history full unstable plaque conditionals with
$\mu_U$ for $\beta$-almost every history in each box. On common refinements with subordinate
partitions, uniqueness of disintegration identifies the Qian--Xie
full and fast plaque conditionals with normalized restrictions of
$\mu_U$ and its vertical conditionals, respectively.
Restrictions and normalizations on the common refinements preserve
almost-everywhere local dimensions. The projection $\pi_0$ is a local
diffeomorphism on each full unstable plaque, so it also preserves local
dimensions. Thus $\mu_U$ is exact-dimensional of dimension $\delta_u$. Its vertical conditionals are
$\mu_z|_{U_z}/\mu_z(U_z)$ for $\zeta_U$-almost every $z$, where
$U_z:=\{w:(z,w)\in U\}$, and have dimension $\delta_f$.
The countable covering by boxes therefore gives the asserted
dimensions of $\mu$ and $\mu_z$.

We now prove the transverse formula. Write $\underline d_\nu$ and
$\overline d_\nu$ for the lower and upper pointwise dimensions of a
measure $\nu$. Fix a positive measure box.
By Fubini's theorem, choose a history for which the preceding
identifications and the Qian--Xie conclusions hold at
$\mu_U$-almost every point.
On countably many small pieces with fixed Qian--Xie chart,
return time and inverse branch, the transverse coordinate of
\cite[Section 6.4, Proposition 6.12]{QX} has the form
$H\circ g_b\circ\pi$.
Here $g_b$ is a local inverse branch of an iterate of $p$, and
$H$ is the transverse change induced by the chosen straightening map. The Qian--Xie coordinate map sends vertical plaques to coordinate
plaques, and both the map and its inverse are Lipschitz.
On smaller product boxes in the straightened coordinates, the induced
transverse change $H$ is therefore locally bi-Lipschitz.
Thus the transverse coordinate is
locally bi-Lipschitz equivalent to $z$.
After a disjoint refinement covering a set of full $\mu_U$-measure,
the base projections of the restrictions of $\mu_U$ to these pieces
are absolutely continuous with respect to $\zeta_U$ and sum to $\zeta_U$. After identifying the transverse coordinates, each projected
restriction is also absolutely continuous with respect to the
corresponding Qian--Xie transverse measure. Lebesgue differentiation therefore identifies the lower and upper local dimensions of each
projected restriction with those of $\zeta_U$ and of the corresponding
Qian--Xie transverse measure almost everywhere with respect to that
restriction.
Applying \cite[Proposition 4.12]{QX} with $i=2$ and letting the
error parameter tend to zero along a countable sequence gives
$\underline d_{\zeta_U}\ge(2L-h_1)/\lambda_b$ almost everywhere
with respect to $\zeta_U$.
By \cite[Lemma 2.2 (i)]{BH}, the exact dimensions $\delta_u$
and $\delta_f$ give
$\overline d_{\zeta_U}\le\delta_u-\delta_f=(2L-h_1)/\lambda_b$
almost everywhere with respect to $\zeta_U$.
Thus $\zeta_U$ is exact-dimensional with the asserted dimension.
\end{proof}

\begin{proof}[Proof of Proposition \ref{prop:fiber-dominant}]
We follow the dimension argument in
\cite[proof of Proposition 5.3]{BH}. Choose a positive measure box from
Lemma \ref{lem:nonuniform-history-boxes}.
Since $\zeta_U=\pi_*\mu_U\ll\mu_p$, the Ma\~n\'e--Manning
formula \cite{Mane1988,Manning1984} and the Lebesgue
differentiation theorem give $\Dim\zeta_U=L/\lambda_b$:
the ratio $\zeta_U(B(z,r))/\mu_p(B(z,r))$ tends to a positive
finite limit for $\zeta_U$-almost every $z$.
Comparing with Proposition \ref{prop:fibered-dimensions} yields
$(2L-h_1)/\lambda_b=L/\lambda_b$, hence $h_1=L$.
Substitution into the formula for $\Dim\mu$ in
Proposition \ref{prop:fibered-dimensions} proves dimension equality.
\end{proof}

\bibliographystyle{amsplain}
\bibliography{base_ref}

\end{document}